\documentclass[a4paper,10pt]{amsart}
\usepackage{amssymb}

\numberwithin{equation}{section}
\newtheorem{theorem}{Theorem}[section]

\newtheorem{lem}[theorem]{Lemma}
\newtheorem{cor}[theorem]{Corollary}

\theoremstyle{remark}
\newtheorem{rem}[theorem]{Remark}

\newcommand{\R}{\mathbb{R}}

\newcommand{\N}{\mathbb{N}}

\author[J.~Benameur]{Jamel Benameur}
\author[M.~Bennaceur]{Mariem Bennaceur}
\address{Higher Institute of Applied Sciences and Technologies of Gab\`es(ISSAT), University of Gab\`es, Tunisia}
\email{\sl jamelbenameur@gmail.com}
\address{Department of Mathematics, Faculty of Science of Gab\`es,
Gab\`es, Tunisia}
\email{\sl Mariemnaceurmariem@outlook.com}

\title[Large  time behaviour of solutions to the 3D-NSE  in $\mathcal X^{\sigma}$ spaces]
{Large  time behaviour of solutions to the 3D-NSE  in $\mathcal X^{\sigma}$ spaces}

\begin{document}
\begin{abstract}
In this paper we study the incompressible Navier-Stokes equations in $L^2(\R^3)\cap\mathcal X^{-1}(\R^3)$.  In the global existence case, we establish that if the solution $u$ is in the space $C(\R^+,L^2\cap\mathcal X^{-1})$, then for $\sigma>-3/2$ the decay of $\|u(t)\|_{\mathcal X^\sigma}$ is at least of the order of $t^{-\frac{\sigma+\frac{3}{2}}{2}}$.  Fourier analysis and standard techniques are used.
\end{abstract}



\maketitle
\tableofcontents


\section{\bf Introduction}

The $3D$ incompressible Navier-Stokes equations are given by:

$$(NS)
  \begin{cases}
     \partial_t u
 -\nu\Delta u+ u.\nabla u   =\;\;-\nabla p\hbox{ in } \mathbb R^+\times \mathbb R^3\\
     {\rm div}\, u = 0 \hbox{ in } \mathbb R^+\times \mathbb R^3\\
    u(0,x) =u^0(x) \;\;\hbox{ in }\mathbb R^3,
  \end{cases}
$$
where $\nu>0$ is the viscosity of fluid, $u=u(t,x)=(u_1,u_2,u_3)$ and $p=p(t,x)$ denote respectively the unknown velocity and the unknown pressure of the fluid at the point $(t,x)\in \mathbb R^+\times \mathbb R^3$, and $(u.\nabla u):=u_1\partial_1 u+u_2\partial_2 u+u_3\partial_3u$, while $u^0=(u_1^o(x),u_2^o(x),u_3^o(x))$ is an initial given velocity. If $u^0$ is quite regular, the divergence free condition determines the pressure $p$.\\
The Navier-Stokes system has the following scaling property : If $u=u(x,t)$ is a solution of $(NS)$ with initial  date $u^0=u^{0}(x)$ on the interval $[0,T]$, then for all $\lambda >0, u_{\lambda}=\lambda u(\lambda^{2}t,\lambda x)$ is a solution of $(NS)$ with initial date $ u_{\lambda}(0,x)=\lambda u^0(\lambda x)$ on the interval $[0,\frac{T}{\lambda^2}].$\\
A fonctional space $(X,\|.\|_X)$ is called critical space of $(NS)$ system if $$\|f_\lambda\|_X=\|f\|_X ;\;\; \forall \lambda>0,\;\; \forall f\in X, $$ where $$ f_\lambda(x)=\lambda f(\lambda x).$$
Particularly, $L^3(\R^3)$, $\dot H^{1/2}(\R^3)$ and $\mathcal {X}^{-1}(\R^3)  $ are critical spaces  for the system $(NS)$. In order to explain the idea of studying the $(NS)$ system in the space $L^2\cap \mathcal X^{-1}$, we introduce the following notation : Two functional spaces $(X,\|.\|_X)$ and $(Y,\|.\|_Y)$ are called "have the same scaling" if, there is a real number $\alpha$ such that
$$\begin{array}{lcl}
\|f_\lambda\|_X&=&\lambda^\alpha\|f\|_X,\;\forall (\lambda,f)\in(0,\infty)\times X\\
\|g_\lambda\|_Y&=&\lambda^\alpha\|g\|_Y,\;\forall (\lambda,g)\in(0,\infty)\times Y.
\end{array}$$
In this case we note $X\approx Y$. For example:
$$\begin{array}{lcl}
\dot H^s(\R^3)&\approx&L^p(\R^3),\;\frac{1}{p}+\frac{s}{3}=\frac{1}{2},\;0\leq s<3/2\;\\
\mathcal X^\sigma(\R^3)&\approx&\dot H^{\sigma+\frac{3}{2}}(\R^3),\;\forall \sigma\in\R.
\end{array}$$
The second is a counter-example of two functional spaces that have the same scaling and are not comparable (see \cite{JB} for $\dot H^{1/2}(\R^3)$ and $\mathcal X^{-1}(\R^3)$). Now, We are ready to give the motivation for this work : Inspired by the works  \cite{ME}, \cite{BS} and \cite{KM} where they proved the decay results of a global solution of $(NS)$ in  homogeneous Sobolev spaces by starting from the $H^s=L^2\cap \dot H^s$ solutions. Here we study the Navier-Stokes system $(NS)$ starting from the $L^2\cap \mathcal X^{-1}$ solutions and proving some optimal decay results. Our first result is the following.
\begin{theorem}\label{theo1} Let $u^0\in\mathcal X^{-1}(\R^3)\cap L^2(\R^3)$ be a divergence free vector fields, then there is a time $T>0$ and unique solution
$u\in C([0,T],\mathcal {X}^{-1}(\R^3)\cap L^2(\R^3))$. Moreover $u\in L^1([0,T],\mathcal X^{1}(\R^3))$. If $\|u^0\|_{\mathcal X^{-1}}<\nu$, then $u$ is global.
\end{theorem}
\begin{rem}\label{rem1} {\bf(i)} If the maximal time $T^*$ is finite then $\int_{0}^{T^*}\|u(t)\|_{\mathcal X^{1}}=+\infty.$ Indeed : The integral form of the system $(NS)$: $$ u(t)=e^{\nu t \Delta}u^0-\int_0^t e^{\nu (t-z)\Delta}\mathbb P(u.\nabla u)dz$$
implies
$$\begin{array}{lcl}
\|u(t)\|_{L^2}&\leq& \displaystyle\|e^{\nu t \Delta}u^0\|_{L^2}+\int_0^t\|e^{\nu (t-z)\Delta}\mathbb P(u.\nabla u)\|_{L^2}dz\\
&\leq&\displaystyle \|u^0\|_{L^2}+\int_0^t\|u\nabla u\|_{L^2}dz\\
&\leq& \displaystyle\|u^0\|_{L^2}+\int_0^t\|u\|_{L^2} \|\nabla u\|_{L^\infty}dz.
\end{array}$$
Using the fact $\|\nabla u\|_{L^\infty}\leq (2\pi)^{-3}\|u\|_{\mathcal X^{1}}$ and  Gronwall lemma we get
\begin{eqnarray}\label{uq}
\|u(t)\|_{L^2}\leq\|u^0\|_{L^2}\exp\big((2\pi)^{-3}\int_0^t \|u\|_{\mathcal X^{1}}\big).
\end{eqnarray}
Then, if $\int_0^{T^*} \|u\|_{\mathcal X^{1}}$ is finite we get $u \in C([0,T^*), L^2 \cap\mathcal X^{-1})\cap L^\infty([0,T^*),L^2 \cap\mathcal X^{-1})$.
Then the solution lives beyond the time $T^*$ which contradicts the fact that $T^*$ is the maximum time of existence.\\
{\bf(ii)} If $\|u^0\|_{\mathcal {X}^{-1}}<\nu/2 $, the above remark and \cite{JB} imply the global existence of solution $u$ of $(NS)$ with $u\in C(\R^+,\mathcal X^{-1})\cap L^1(\R^+,\mathcal X^{1})\cap C(\R^+,L^2) $. Moreover,
\begin{eqnarray}\label{eQ2}
\|u(t)\|_{\mathcal {X}^{-1}}+\frac{\nu}{2} \int_{0}^{t}\| u\|_{\mathcal {X}^{1}}\leq \|u^0\|_{\mathcal {X}^{-1}};\;\; \forall t\geq 0.
\end{eqnarray}
{\bf(iii)} Using (i)-(ii) and \cite{JB}, we get if $u \in C(\R^+, L^2 \cap\mathcal X^{-1})$ is a global solution of $(NS)$, then $u\in L^1(\R^+,\mathcal X^1(\R^3))$.\\
{\bf(iv)} Using (i)-(ii)-(iii) and \cite{JB}, we get if $u \in C(\R^+, L^2 \cap\mathcal X^{-1})$ is a global solution of $(NS)$, then $u\in \mathcal C_b(\R^+,L^2(\R^3))$. Indeed: By \cite{JB} there is a time $t_0\geq0$ such that $\|u(t_0)\|_{\mathcal X^{-1}}<\nu/2$. Then (i)-(ii) imply for $t\geq t_0$
$$\begin{array}{lcl}
\|u(t)\|_{L^2}&\leq& \displaystyle\|u(t_0)\|_{L^2}\exp\big((2\pi)^{-3}\int_{t_0}^\infty\|u(z)\|_{\mathcal X^1}dz\big)\\
&\leq& \displaystyle\|u(t_0)\|_{L^2}\exp\big((2\pi)^{-3}\frac{2\|u(t_0)\|_{\mathcal X^{-1}}}{\nu}\big)\\
&\leq& \displaystyle\|u(t_0)\|_{L^2}\exp\big((2\pi)^{-3}\big)\\
&\leq& 2\|u(t_0)\|_{L^2},
\end{array}$$
which implies
$$\|u(t)\|\leq2\max_{0\leq z\leq t_0}\|u(z)\|_{L^2},\;\forall t\geq0.$$
Particularly, if $\|u^0\|_{\mathcal X^{-1}}<\nu/2$ we get
$$\|u(t)\|_{L^2}\leq2\|u^0\|_{L^2},\;\forall t\geq0.$$
\end{rem}
Before stating the result of decay of the global solution of $(NS)$, we recall the following results which will be useful in the following
\begin{theorem}\label{th13}(see \cite{JB})
Let $u\in C(\R^+,\mathcal X^{-1}(\R^3))$ be a global solution of Navier-Stokes system. Then $$\lim_{t\rightarrow \infty}\|u(t)\|_{\mathcal X^{-1}}=0.$$
\end{theorem}
\begin{theorem}\label{th14}(see \cite{HB}) There exists a positive constant $\epsilon_{0}> 0$ such that for any initial data $u^0$ in $ \mathcal X^{-1}(\R^3)$ with $\|u^0\|_{\mathcal X^{-1}}<\epsilon_{0}$, the solution  of Navier-Stokes system is analytic in the sense that
$$\|\exp(\sqrt{\nu t}|D|)u(t)\|_{\mathcal X^{-1}}+ \frac{\nu}{2}\int^t_0\|\exp(\sqrt{\nu z}|D|)u(z) \|_{\mathcal X^{1}}dz\leq 2\|u^0\|_{\mathcal X^{-1}},\;\forall t\geq0.$$
\end{theorem}
\begin{theorem}\label{th15}(see \cite{TK})
For any initial data $u^0\in H^s(\R^3)$ with $\rm {div}\,u^0=0$, there exists a unique solution $u\in C([0,T_0],H^s(\R^3))$ such that $T_{0}=T_0(s,\|u^0\|_{H^s}).$
\end{theorem}
Our second result is the following.
\begin{theorem} \label{theo2}
Let $u\in{\mathcal C}(\R^+,\mathcal X^{-1}(\R^3)\cap L^2(\R^3))$ be a global solution of Navier-Stokes system. Then
\begin{equation}\label{th2eq1}
\lim_{t\rightarrow +\infty}\|u(t)\|_{\mathcal X^{-1}\cap L^2}=0,
\end{equation}
Precisely,
\begin{equation}\label{th2eq2}
\|u(t)\|_{\mathcal X^{-1}} = o(t^{-\frac{1}{4}}) ; \;\; t\rightarrow +\infty .
\end{equation}
\end{theorem}
Using theorem \ref{theo2} and theorem \ref{th14} which characterizes the regularizing effect of the Navier-Stokes equations, we get the following decay result of $\|u(t)\|_{\mathcal X^\sigma}$.
\begin{cor}\label{cor1}
Let  $u\in{\mathcal C}(\R^+,\mathcal X^{-1}(\R^3)\cap L^2(\R^3))$ be a global solution of Navier-Stokes system. Then, for all $\sigma>-3/2$, we have $u \in {\mathcal C}((0,+\infty),\mathcal X^{\sigma})$ and
\begin{equation*}
\|u(t)\|_{\mathcal X^{\sigma}}=o(t^{-\frac{\sigma+3/2}{2}}) ;\;\; t\rightarrow +\infty.
\end{equation*}
\end{cor}

The remainder of our paper is organized as follows. In the second section we give some notations, definitions and preliminary results.  Section 3 is devoted to prove the well posedness of $(NS)$ in $L^2\cap \mathcal X^{-1}$ space, this proof used the Fixed Point Theorem with a good choice of space $X=\mathcal C([0,T],L^2\cap\mathcal X^{-1})\cap L^1([0,T],\mathcal X^1)$. In section 4 we prove the decay of global solutions in $L^2\cap \mathcal X^{-1}$, this proof used a Fourier analysis and standard techniques. Section 5 is devoted to prove the decay results of the global solution in $\mathcal X^\sigma$, this proof uses in a fundamental way the decay in $L^2\cap \mathcal X^{-1}$.
\section{\bf Notations and preliminary results}
\subsection{Notations} In this section, we collect some notations and definitions that will be used later.\\
$\bullet$ The Fourier transformation is normalized as
$$
\mathcal{F}(f)(\xi)=\widehat{f}(\xi)=\int_{\mathbb R^3}\exp(-ix.\xi)f(x)dx,\,\,\,\xi=(\xi_1,\xi_2,\xi_3)\in\mathbb R^3.
$$
$\bullet$ The inverse Fourier formula is
$$
\mathcal{F}^{-1}(g)(x)=(2\pi)^{-3}\int_{\mathbb R^3}\exp(i\xi.x)g(\xi)d\xi,\,\,\,x=(x_1,x_2,x_3)\in\mathbb R^3.
$$
$\bullet$ The convolution product of a suitable pair of function $f$ and $g$ on $\mathbb R^3$ is given by
$$
(f\ast g)(x):=\int_{\mathbb R^3}f(y)g(x-y)dy.
$$
$\bullet$ If $f=(f_1,f_2,f_3)$ and $g=(g_1,g_2,g_3)$ are two vector fields, we set
$$
f\otimes g:=(g_1f,g_2f,g_3f),
$$
and
$$
{\rm div}\,(f\otimes g):=({\rm div}\,(g_1f),{\rm div}\,(g_2f),{\rm div}\,(g_3f)).
$$
Moreover, if $\rm{div}\,g=0$ we obtain
$$
{\rm div}\,(f\otimes g):=g_1\partial_1f+g_2\partial_2f+g_3\partial_3f:=g.\nabla f.
$$
$\bullet$ Let $(B,||.||)$, be a Banach space, $1\leq p \leq\infty$ and  $T>0$. We define $L^p_T(B)$ the space of all
measurable functions $[0,t]\ni t\mapsto f(t) \in B$ such that $t\mapsto||f(t)||\in L^p([0,T])$.\\
$\bullet$ The Sobolev space $H^s(\R^3)=\{f\in \mathcal S'(\R^3);\;(1+|\xi|^2)^{s/2}\widehat{f}\in L^2(\R^3)\}$.\\
$\bullet$ The homogeneous Sobolev space $\dot H^s(\R^3)=\{f\in \mathcal S'(\R^3);\;\widehat{f}\in L^1_{loc}\;{\rm and}\;|\xi|^s\widehat{f}\in L^2(\R^3)\}$.\\
$\bullet$ The Lei-Lin space $\mathcal X^\sigma(\R^3)=\{f\in \mathcal S'(\R^3);\;\widehat{f}\in L^1_{loc}\;{\rm and}\;|\xi|^\sigma\widehat{f}\in L^1(\R^3)\}$.\\
\subsection{Preliminary results}
In this section, we recall some classical results and we give new technical lemmas.\\
\begin{lem}\label{lem1}
We have $\mathcal X^{-1}(\R^3)\cap \mathcal X^{1}(\R^3)\hookrightarrow \mathcal X^0(\R^3)$. Precisely, we have
\begin{eqnarray}\label{enq1}
\|f\|_{\mathcal{X}^0(\mathbb{R}^3)}\leq\|f\|_{\mathcal{X}^{-1}(\mathbb{R}^3)}^{1/2}\|f\|_{\mathcal{X}^{1}(\mathbb{R}^3)}^{1/2},\;\forall f\in\mathcal X^{-1}(\R^3)\cap \mathcal X^{1}(\R^3).
\end{eqnarray}
\end{lem}
\noindent{\bf Proof.}
We can write
$$\begin{array}{lcl}
\|f\|_{\mathcal X^0}&=&\displaystyle\int_{\R^3} |\widehat{f}(\xi)|d\xi\\
              &\leq&\displaystyle \int_{\R^3} |\xi|^{-1/2}|\widehat{f}(\xi)|^{1/2}\frac{|\widehat{f}(\xi)|^{1/2}}{|\xi|^{1/2}}d\xi.
\end{array}$$
Cauchy-Schwartz inequality gives the result.
\begin{lem}\label{lem2}
Let $\sigma,s\in\R$ such that $0<\sigma+\frac{3}{2}<s$. Then $H^s(\R^3)\hookrightarrow \mathcal X^\sigma(\R^3)$. Precisely, there is a constant $C=C(s,\sigma)$ such that
\begin{eqnarray}\label{enq1}
\|f\|_{\mathcal{X}^\sigma(\mathbb{R}^3)}\leq C
\|f\|_{L^2(\mathbb{R}^3)}^{1-\frac{\sigma+\frac{3}{2}}{s}}\|f\|_{\dot H^s(\mathbb{R}^3)}^{\frac{\sigma+\frac{3}{2}}{s}},\; \;\forall f\in H^s(\R^3).
\end{eqnarray}
\end{lem}
\noindent{\bf Proof.} For $\lambda>0$, we have
$$\|f\|_{\mathcal X^\sigma}=I_\lambda+J_\lambda,$$
with
$$\begin{array}{lcl}
I_\lambda&=&\int_{|\xi|<\lambda}|\xi|^\sigma|\widehat{f}(\xi)|d\xi\\
J_\lambda&=&\int_{|\xi|>\lambda}|\xi|^\sigma|\widehat{f}(\xi)|d\xi.
\end{array}$$
We have
$$\begin{array}{lcl}
I_\lambda&\leq&\Big(\int_{|\xi|<\lambda}|\xi|^{2\sigma}d\xi\Big)^{1/2}\|f\|_{L^2}\\
&\leq&\frac{c}{\sqrt{2\sigma+3}}\lambda^{\sigma+\frac{3}{2}}\|f\|_{L^2}\end{array}$$
$$\begin{array}{lcl}
J_\lambda&\leq&\Big(\int_{|\xi|>\lambda}|\xi|^{2(\sigma-s)}d\xi\Big)^{1/2}\|f\|_{\dot H^s}\\
&\leq&\Big(\int_{|\xi|>\lambda}|\xi|^{2(\sigma-s)}d\xi\Big)^{1/2}\|f\|_{\dot H^s}\\
&\leq&\frac{C}{\sqrt{s-\sigma-\frac{3}{2}}}\lambda^{\sigma+\frac{3}{2}-s}\|f\|_{\dot H^s}.
\end{array}$$
For $\lambda=(\|f\|_{\dot H^s}/\|f\|_{L^2})^{1/s}$, we obtain the desired result.
\begin{lem}
Let $\sigma_{0}>-3/2 .$
If we have $$\mathcal X^{\sigma_{0}}(\R^3) \cap  L^2(\R^3)\hookrightarrow \mathcal X^{\sigma}(\R^3); \;\;\;\forall -3/2<\sigma\leq\sigma_{0} .$$
Precisely
\begin{eqnarray}
\|f\|_{\mathcal{X}^\sigma}\leq c_{0}
\|f\|_{L^2}^{1-\theta}\|f\|_{\mathcal{X}^{\sigma_{0}}}^{\theta},\; \;\;\forall   c_{0}=c(\sigma_{0},\sigma),\;\; \theta=\frac{\sigma+\frac{3}{2}}{\frac{3}{2}+\sigma_{0}}.
\end{eqnarray}
\end{lem}
\noindent{\bf Proof.} For $\lambda>0$, we have
$$\|f\|_{\mathcal{X}^\sigma}=A(\lambda)+B(\lambda),$$
with
$$\begin{array}{lcl}
A(\lambda)&=&\int_{\{\xi\in\R^3/|\xi|<\lambda\}}|\xi|^\sigma|\widehat{f}(\xi)|d\xi\\
B(\lambda)&=&\int_{\{\xi\in\R^3/|\xi|>\lambda\}}|\xi|^\sigma|\widehat{f}(\xi)|d\xi.
\end{array}$$
We have
$$\begin{array}{lcl}
A(\lambda)&\leq&\displaystyle\int_{\{\xi\in\R^3/|\xi|<\lambda\}}|\xi|^{2\sigma}d\xi\Big)^{1/2}\|f\|_{L^2}\\
&\leq&\displaystyle\frac{c}{\sqrt{2\sigma+3}}\lambda^{\sigma+\frac{3}{2}}\|f\|_{L^2}\\\\
B(\lambda)&\leq&\displaystyle\int_{\{\xi\in\R^3/|\xi|>\lambda\}}|\xi|^{(\sigma-\sigma_{0})}|\xi|^{\sigma_{0}}|\widehat{f}(\xi)|d\xi\\
&\leq&\displaystyle\lambda^{\sigma-\sigma_{0}}\|f\|_{\mathcal{X}^{\sigma_{0}}},
\end{array}$$
which imply
$$\|f\|_{\mathcal{X}^\sigma} \leq \frac{c}{\sqrt{2\sigma+3}}\lambda^{\sigma+\frac{3}{2}}\|f\|_{L^2}+\lambda^{\sigma-\sigma_{0}}\|f\|_{\mathcal{X}^{\sigma_{0}}}.$$
For $\lambda=(\|f\|_{{\mathcal{X}^{\sigma_{0}}}}/\|f\|_{L^2})^{1/(3/2+\sigma_{0})}$,
 we obtain
$$\|f\|_{\mathcal{X}^\sigma}\leq c_{\sigma,\sigma_{0}}
\|f\|_{L^2}^{\frac{\sigma_{0}-\sigma}{\frac{3}{2}+\sigma_{0}}}\|f\|_{\mathcal{X}^{\sigma_{0}}}^{\frac{\sigma+\frac{3}{2}}{\frac{3}{2}+\sigma_{0}}}.$$
\begin{lem}\label{lem3}
Let $f,g\in L^\infty_T(\mathcal X^{-1}(\R^3)\cap L^{2}(\R^3))\cap L^1_T(\mathcal X^{1}(\R^3))$ such that ${\rm{div}}\,f=0$ almost everywhere. Then
\begin{eqnarray}\label{enq1}
\sup_{0\leq t\leq T}\|\int_0^te^{\nu (t-z)\Delta}\mathbb P(f.\nabla g)dz\|_{\mathcal{X}^{-1}}\leq\|f\|_{L^\infty_T(\mathcal{X}^{-1})}^{1/2}\|f\|_{L^1_T(\mathcal{X}^{1})}^{1/2}
\|g\|_{L^\infty_T(\mathcal{X}^{-1})}^{1/2}\|g\|_{L^1_T(\mathcal{X}^{1})}^{1/2},
\end{eqnarray}
\begin{eqnarray}\label{enq2}
\sup_{0\leq t\leq T}\|\int_0^te^{\nu (t-z)\Delta}\mathbb P(f.\nabla g)dz\|_{L^2}\leq(2\pi)^{-3}\|f\|_{L^\infty_T(L^2)}\|g\|_{L^1_T(\mathcal{X}^{1})},
\end{eqnarray}
\begin{eqnarray}\label{enq3}
\int_0^T\|\int_0^te^{\nu(t-z)\Delta}\mathbb P(f.\nabla g)dz\|_{\mathcal{X}^1}dt\leq\nu^{-1}\|f\|_{L^\infty_T(\mathcal{X}^{-1})}^{1/2}\|f\|_{L^1_T(\mathcal{X}^{1})}^{1/2}
\|g\|_{L^\infty_T(\mathcal{X}^{-1})}^{1/2}\|g\|_{L^1_T(\mathcal{X}^{1})}^{1/2}.
\end{eqnarray}
\end{lem}
\noindent{\bf Proof.}\\
$\bullet$ Proof of (\ref{enq1}): We can write
$$\begin{array}{lcl}
\|\int_0^te^{\nu (t-z)\Delta}\mathbb P(f.\nabla g)dz\|_{\mathcal X^{-1}}&\leq&\int_0^t\|e^{\nu (t-z)\Delta}\mathbb P(f.\nabla g)\|_{\mathcal X^{-1}}dz\\
&\leq&\int_0^t\|f.\nabla g\|_{\mathcal X^{-1}}dz\\
&\leq&\int_0^T\|{\rm{div}}\,(f\otimes g)\|_{\mathcal X^{-1}}dz\\
&\leq&\int_0^T\|f\otimes g\|_{\mathcal X^{0}}dz\\
&\leq&\int_0^T\|f\|_{\mathcal X^{0}}\|g\|_{\mathcal X^{0}}dz\\
&\leq&\int_0^T\|f\|_{\mathcal X^{-1}}^{1/2}\|f\|_{\mathcal X^{1}}^{1/2}\|g\|_{\mathcal X^{-1}}^{1/2}\|g\|_{\mathcal X^{1}}^{1/2}dz\\
&\leq&\|f\|_{L^\infty_T(\mathcal X^{-1})}^{1/2}\|g\|_{L^\infty_T(\mathcal X^{-1})}^{1/2}\int_0^T\|f\|_{\mathcal X^{1}}^{1/2}\|g\|_{\mathcal X^{1}}^{1/2}dz\\
&\leq&\|f\|_{L^\infty_T(\mathcal{X}^{-1})}^{1/2}\|f\|_{L^1_T(\mathcal{X}^{1})}^{1/2}
\|g\|_{L^\infty_T(\mathcal{X}^{-1})}^{1/2}\|g\|_{L^1_T(\mathcal{X}^{1})}^{1/2}.
\end{array}$$
\noindent$\bullet$  Proof of (\ref{enq2}): We can write
$$\begin{array}{lcl}
\|\int_0^te^{\nu (t-z)\Delta}\mathbb P(f.\nabla g)dz\|_{L^2}&\leq&\int_0^t\|e^{\nu (t-z)\Delta}\mathbb P(f.\nabla g)\|_{L^2}dz\\
&\leq&\int_0^t\|f.\nabla g\|_{L^2}dz\\
&\leq&\int_0^T\|f.\nabla g\|_{L^2}dz\\
&\leq&\int_0^T\|f\|_{L^2}\|\nabla g\|_{L^\infty}dz\\
&\leq&(2\pi)^{-3}\|f\|_{L^\infty_T(L^2)}\int_0^T\|g\|_{\mathcal X^1}dz\\
&\leq&(2\pi)^{-3}\|f\|_{L^\infty_T(L^2)}\|g\|_{L^1_T(\mathcal X^1)}dz.
\end{array}$$
\noindent$\bullet$  Proof of (\ref{enq3}): We can write
$$\begin{array}{lcl}
\int_0^T\|\int_0^te^{\nu (t-z)\Delta}\mathbb P(f.\nabla g)dz\|_{\mathcal X^{1}}dt&\leq&\int_0^T\int_0^t\|e^{\nu (t-z)\Delta}\mathbb P(f.\nabla g)\|_{\mathcal X^{1}}dzdt\\
&\leq&\int_0^T\int_0^t\int_{\R^3}e^{-\nu (t-z)|\xi|^2}|\xi|.|\mathcal F({\rm{div}}\,(f\otimes g))(z,\xi)|d\xi dzdt\\
&\leq&\int_0^T\int_0^t\int_{\R^3}e^{-\nu (t-z)|\xi|^2}|\xi|^2.|\mathcal F(f\otimes g)(z,\xi)|d\xi dzdt\\
&\leq&\int_{\R^3}|\xi|^2\Big(\int_0^T\int_0^te^{-\nu (t-z)|\xi|^2}|\mathcal F(f\otimes g)(z,\xi)|dzdt\Big)d\xi\\
&\leq&\int_{\R^3}|\xi|^2\Big(\|e^{-\nu t|\xi|^2}*_t|\mathcal F(f\otimes g)(t,\xi)\|_{L^1([0,T])}\Big)d\xi\\
&\leq&\int_{\R^3}|\xi|^2\Big(\|e^{-\nu .|\xi|^2}\|_{L^1([0,T])}\|\mathcal F(f\otimes g)(.,\xi)\|_{L^1([0,T])}\Big)d\xi\\
&\leq&\int_{\R^3}|\xi|^2\Big(\frac{1-e^{-\nu T|\xi|^2}}{\nu|\xi|^2}\int_0^T|\mathcal F(f\otimes g)(t,\xi)|dt\Big)d\xi\\
&\leq&\nu^{-1}\int_0^T\int_{\R^3}|\mathcal F(f\otimes g)(t,\xi)|d\xi dt\\
&\leq&\nu^{-1}\int_0^T\|f\otimes g(t)\|_{\mathcal X^0}dt\\
&\leq&\nu^{-1}\int_0^T\|f\otimes g\|_{\mathcal X^{0}}dz\\
&\leq&\nu^{-1}\int_0^T\|f\|_{\mathcal X^{0}}\|g\|_{\mathcal X^{0}}dz\\
&\leq&\nu^{-1}\int_0^T\|f\|_{\mathcal X^{-1}}^{1/2}\|f\|_{\mathcal X^{1}}^{1/2}\|g\|_{\mathcal X^{-1}}^{1/2}\|g\|_{\mathcal X^{1}}^{1/2}dz\\
&\leq&\nu^{-1}\|f\|_{L^\infty_T(\mathcal X^{-1})}^{1/2}\|g\|_{L^\infty_T(\mathcal X^{-1})}^{1/2}\int_0^T\|f\|_{\mathcal X^{1}}^{1/2}\|g\|_{\mathcal X^{1}}^{1/2}dz\\
&\leq&\nu^{-1}\|f\|_{L^\infty_T(\mathcal{X}^{-1})}^{1/2}\|f\|_{L^1_T(\mathcal{X}^{1})}^{1/2}
\|g\|_{L^\infty_T(\mathcal{X}^{-1})}^{1/2}\|g\|_{L^1_T(\mathcal{X}^{1})}^{1/2}.
\end{array}$$
\begin{lem}\label{lemt1}
Let $T>0$ and $f :  [0,T]\rightarrow \R_{+} $ be continuous function such that
\begin{equation}\label{eql1}
f(t)\leq M_{0}+\theta_{1}f(\theta_{2}t)  ;  \;\;\; \forall \; 0 \leq  t \leq T .
\end{equation}
with $M_{0}\geq 0$ and $\theta_{1},\theta_{2} \in (0,1).$ Then
\begin{equation*}
f(t)\leq \frac{M_{0}}{1-\theta_{1}}  ; \;\;\;  \forall  0 \leq t\leq T.
\end{equation*}
\end{lem}
\noindent{\bf Proof.}
As $f$ is a positive and continuous function, then there  is a time $t_{0}\in [0,T]$ such that $$ 0 \leq f(t_{0})= \max_{0\leq t\leq T}f(t).$$
Applying (\ref{eql1}) at $ t=t_{0}$ we get
$$f(t_{0})\leq M_{0}+\theta_{1}f(\theta_{2}t_{0})\leq  M_{0}+\theta_{1}f(t_{0})$$
which implies $f(t_{0}) \leq\frac{ M_{0}}{1-\theta_{1}}.$ As $ f(t_{0}) = \max_{0\leq t\leq T}f(t )$, we get the desired result.
\begin{rem}\label{remlemt1} Applying Lemma \ref{lemt1} to a positive continuous function $f: \R_{+}\rightarrow \R_{+}$ satisfying
$$
f(t)\leq M_{0}+\theta_{1}f(\theta_{2}t) ; \;\;\; \forall  t\geq0
$$
with $ M_{0}\geq 0$ and $\theta_{1},\theta_{2} \in (0,1)$, we obtain
\begin{equation*}
 \limsup_{t\rightarrow +\infty} f(t) \leq \frac{M_{0}}{1-\theta_{1}}.
\end{equation*}
\end{rem}
\section{\bf Well posedness results in $L^2(\R^3)\cap \mathcal X^{-1}(\R^3)$} In this section we prove Theorem \ref{theo1}.
To prove the existence result we need the following remark : For $f\in L^2(\R^3)\cap \mathcal X^{-1}(\R^3)$ and $\varepsilon_0>0$ there is $\lambda>0$ such that
$$\|\lambda f(\lambda .)\|_{\mathcal X^{-1}}=\|f\|_{\mathcal X^{-1}}\;\;{\rm{and}}\;\;\|\lambda f(\lambda .)\|_{L^2}<\varepsilon_0.$$
Precisely, just take $\lambda=\frac{\varepsilon_0^2}{4\|f\|_{L^2}^2+1}$. Then we can choose $\lambda_0>0$ such that
$$\|\lambda_0 u^0(\lambda_0 .)\|_{\mathcal X^{-1}}=\|u^0\|_{\mathcal X^{-1}}\;\;{\rm{and}}\;\;\|\lambda_0u^0(\lambda_0 .)\|_{L^2}<\frac{1}{48}.$$
Consider then the Navier-Stokes system
$$(NS_{\lambda_0})
  \begin{cases}
     \partial_t v
 -\nu\Delta v+ v.\nabla v   =\;\;-\nabla q\hbox{ in } \mathbb R^+\times \mathbb R^3\\
     {\rm div}\, v = 0 \hbox{ in } \mathbb R^+\times \mathbb R^3\\
    v(0,x) =\lambda_0u^0(\lambda_0 x) \;\;\hbox{ in }\mathbb R^3.
  \end{cases}
$$
If the system $(NS_{\lambda_0})$ has a unique solution $v$ in $C([0,T],L^2\cap \mathcal X^{-1})$, then $u=\lambda_0^{-1}v(\lambda_0^{-2}t,\lambda_0^{-1      }x)$ is a solution of Navier-Stokes system starting by $u^0$. Therefore, we can assume in the following that
\begin{equation}\label{L2}
\|u^0\|_{L^2}<\frac{1}{48}.
\end{equation}
Let's go back to the proof of Theorem \ref{theo1}. A uniqueness in $L^2\cap \mathcal X^{-1}$ is given by the uniqueness in $\mathcal X^{-1}$,(see \cite{LL}). It remains a proven existence, for this let $k_0\in\N^*$ such that $$\int_{\{\xi\in\R^3/|\xi|>k_0\}}\frac{|\widehat{u^0}(\xi)|}{|\xi|}d\xi<\min(\frac{\nu}{16},\frac{1}{16}).$$
Put
$$\begin{array}{lcl}
a^0&=&\mathcal F^{-1}({\bf 1}_{|\xi|<k}\widehat{u^0}(\xi))\\
b^0&=&\mathcal F^{-1}({\bf 1}_{|\xi|\geq k}\widehat{u^0}(\xi)).
\end{array}$$
We have $a^0\in H^s(\R^3)$,(for all $s\geq0$) and
\begin{equation}\label{eq41}
a^0\in H^s(\R^3),\;\;\forall s\geq0,
\end{equation}
\begin{equation}\label{eq42}
\|b^0\|_{\mathcal X^{-1}}<\min(\frac{\nu}{16},\frac{1}{16}).
\end{equation}
Moreover
\begin{equation}\label{eq412}
\|a^0\|_{L^2}\leq\|u^0\|_{L^2}\;\;{\rm{and}}\;\;\|b^0\|_{L^2}\leq\|u^0\|_{L^2}.
\end{equation}
There is a time $T_0>0$ such that the system $(NS)$ has a unique solution $a$ in $C([0,T_0],H^4(\R^3))$ with initial condition $a^0$ (see \cite{TK}). Using the fact(see Lemma \ref{lem2})
$$H^4(\R^3)\hookrightarrow L^2(\R^3)\cap\mathcal X^{-1}(\R^3)\cap \mathcal X^1(\R^3),$$
we get
\begin{equation}\label{eq43}
a\in C([0,T_0],L^2(\R^3)\cap\mathcal X^{-1}(\R^3)\cap \mathcal X^1(\R^3)).
\end{equation}
Using the regularity of the function $a$ and inequality (\ref{eq412}), we obtain
\begin{equation}\label{eq431}
\|a(t)\|_{L^2}^2+2\nu\int_0^t\|\nabla a(z)\|_{L^2}^2=\|a^0\|_{L^2}^2\leq \|u^0\|_{L^2}^2,\;\forall t\in[0,T_0].
\end{equation}
Put $b=u-a$, $b$ satisfies the following system
$$(RNS)
  \begin{cases}
     \partial_t b
 -\nu\Delta b+ b.\nabla b+b.\nabla a+a.\nabla b   =\;\;-\nabla q\hbox{ in } \mathbb R^+\times \mathbb R^3\\
     {\rm div}\, b = 0 \hbox{ in } \mathbb R^+\times \mathbb R^3\\
    b(0,x) =b^0(x) \;\;\hbox{ in }\mathbb R^3.
  \end{cases}
$$
The integral form of $(RNS)$ is
$$b=\psi(b)=e^{\nu t\Delta}b^0-\int_0^te^{\nu(t-\tau)\Delta}\mathbb P(a.\nabla b)-\int_0^te^{\nu(t-\tau)\Delta}\mathbb P(b.\nabla a)-\int_0^te^{\nu(t-\tau)\Delta}\mathbb P(b.\nabla b).$$
Put
$$\begin{array}{lcl}
f_0&=&e^{\nu t\Delta}b^0\\
L(b)&=&-\int_0^te^{\nu(t-\tau)\Delta}\mathbb P(a.\nabla b)-\int_0^te^{\nu(t-\tau)\Delta}\mathbb P(b.\nabla a)\\
Q(b)&=&-\int_0^te^{\nu(t-\tau)\Delta}\mathbb P(b.\nabla b).
\end{array}$$
For $T>0$ put the space $$X_T=C([0,T],L^2(\R^3)\cap \mathcal X^{-1}(\R^3))\cap L^1([0,T],\mathcal X^1(\R^3)).$$
This vector space is equipped with the norm
$$\|f\|_{\varepsilon,T}=\|f\|_{L^\infty_T(L^2)}+\|f\|_{L^\infty_T(\mathcal X^{-1})}+\|f\|_{L^1_T(\mathcal X^{1})}.$$
For $\varepsilon,T>0$ (to fixed later), such that $T\leq T_0$, put the closed subset of $X_T$ defined by
$$B(\varepsilon,T)=\Big\{f\in X_T;\left\{\begin{array}{lcl}\|f\|_{L^\infty_T(L^2)}&\leq& 2\|b^0\|_{L^2}\\
\|f\|_{L^\infty_T(\mathcal X^{-1})}&\leq& 2\|b^0\|_{\mathcal X^{-1}}\\
\|f\|_{L^1_T(\mathcal X^{1})}&\leq& \varepsilon\end{array}\right.\Big\}$$
Explanation of the choice of $\varepsilon$ and $T$ : We have
$$\begin{array}{lcl}
\|f_0\|_{L^\infty_T(\mathcal X^{-1})}&\leq& \|b^0\|_{\mathcal X^{-1}}\\
\|f_0\|_{L^\infty_T(L^2)}&\leq& \|b^0\|_{L^2}\\
\|f_0\|_{L^1_T(\mathcal X^{1})}&=& \int_0^T\int_{\R^3}e^{-\nu t|\xi|^2}|\xi|.|\widehat{b^0}(\xi)|d\xi dt\\
&=& \int_{\R^3}\Big(\int_0^Te^{-\nu t|\xi|^2}dt\Big)|\xi|.|\widehat{b^0}(\xi)|d\xi\\
&=& \int_{\R^3}\Big(\frac{1-e^{-\nu T|\xi|^2}}{\nu|\xi|^2}dt\Big)|\xi|.|\widehat{b^0}(\xi)|d\xi\\
&=& \int_{\R^3}(1-e^{-\nu T|\xi|^2})\frac{|\widehat{b^0}(\xi)|}{|\xi|}d\xi.
\end{array}$$
Dominate Convergence Theorem implies
\begin{equation}\label{eq44}\lim_{t\rightarrow 0^+}\|f_0\|_{L^1_T(\mathcal X^{1})}=0.\end{equation}
Let $0<\varepsilon<1/24$ and $0<T\leq T_0$ such that
$$\begin{array}{lcl}
&(C1)&\|a\|_{L^\infty_T(\mathcal X^{-1})}^{1/2}\|a\|_{L^1_T(\mathcal X^{1})}^{1/2}\sqrt{2}\sqrt{\varepsilon}\|b^0\|_{\mathcal X^{-1}}^{1/2}
\leq \frac{\|b^0\|_{\mathcal X^{-1}}}{4}\\
&(C2)&\|a^0\|_{L^2}\varepsilon+2\|a\|_{L^1_T(\mathcal X^{1})}\|b^0\|_{L^2}
\leq \frac{\|b^0\|_{L^2}}{4}\\
&(C3)&\nu^{-1}\|a\|_{L^\infty_T(\mathcal X^{-1})}^{1/2}\|a\|_{L^1_T(\mathcal X^{1})}^{1/2}\sqrt{2}\sqrt{\varepsilon}\|b^0\|_{\mathcal X^{-1}}^{1/2}\leq\varepsilon/3\\
&(C4)& \varepsilon+2\|b^0\|_{L^2}\leq 1/12\\
&(C5)& (1+\nu^{-1})2\sqrt{2}\sqrt{\varepsilon}\|b^0\|_{\mathcal X^{-1}}^{1/2}\leq 1/12\\
&(C6)& \|f_0\|_{L^1_T(\mathcal X^{1})}\leq \varepsilon/3\\
&(C7)&2\sqrt{2\varepsilon}\|b_0\|_{\mathcal X^{-1}}\leq 1/12\\
&(C8)&\|a\|_{L^\infty_T(\mathcal X^{-1})}\|a\|_{L^1_T(\mathcal X^{1})}\leq 1/12\\
&(C9)&\|a\|_{L^1_T(\mathcal X^{1})}\leq 1/24.
\end{array}$$
These choices are possible just use the equations (\ref{eq42})-(\ref{eq43})-(\ref{eq44}). Now we want to prepare to apply the Fixed Point Theorem, for this we prove the following
\begin{equation}\label{eq45}\psi(B(\varepsilon,T))\subset B(\varepsilon,T).\end{equation}
\begin{equation}\label{eq46}\|\psi(\alpha_1)-\psi(\alpha_2)\|_{\varepsilon,T}\leq \frac{1}{2}\|\alpha_1-\alpha_2\|_{\varepsilon,T},\;\forall \alpha_1,\alpha_2\in B(\varepsilon,T).\end{equation}
\noindent\underline{\bf Proof of (\ref{eq45}):} Using inequality (\ref{enq1}), we obtain
$$\begin{array}{lcl}
\|L(b)\|_{L^\infty_T(\mathcal X^{-1})}&\leq&\|a\|_{L^\infty_T(\mathcal X^{-1})}^{1/2}\|a\|_{L^1_T(\mathcal X^{1})}^{1/2}\|b\|_{L^\infty_T(\mathcal X^{-1})}^{1/2}\|b\|_{L^1_T(\mathcal X^{1})}^{1/2}\\
&\leq&\|a\|_{L^\infty_T(\mathcal X^{-1})}^{1/2}\|a\|_{L^1_T(\mathcal X^{1})}^{1/2}\sqrt{2}\sqrt{\varepsilon}\|b^0\|_{\mathcal X^{-1}}^{1/2}\\
&\leq&\frac{\|b^0\|_{\mathcal X^{-1}}}{4},\;(by\,(C1))\\\\
\|Q(b)\|_{L^\infty_T(\mathcal X^{-1})}&\leq&\|b\|_{L^\infty_T(\mathcal X^{-1})}\|b\|_{L^1_T(\mathcal X^{1})}\\
&\leq&2\varepsilon\|b^0\|_{\mathcal X^{-1}}\\
&\leq&\frac{\|b^0\|_{\mathcal X^{-1}}}{4}.
\end{array}$$
Then
\begin{equation}\label{eq301}
\|\psi(b)\|_{L^\infty_T(\mathcal X^{-1})}\leq 2\|b^0\|_{\mathcal X^{-1}},\;\forall b\in B(\varepsilon,T).
\end{equation}
Similarly, inequality (\ref{enq2}) gives
$$\begin{array}{lcl}
\|L(b)\|_{L^\infty_T(L^2)}&\leq&\|a\|_{L^\infty_T(L^2)}\|b\|_{L^1_T(\mathcal X^{1})}+\|b\|_{L^\infty_T(L^2)}\|a\|_{L^1_T(\mathcal X^{1})}\\
&\leq&\|a^0\|_{L^2} \varepsilon+2\|a\|_{L^1_T(\mathcal X^{1})}\|b^0\|_{L^2}\\
&\leq&\frac{\|b^0\|_{L^2}}{4},\;(by\,(C2))\\\\
\|Q(b)\|_{L^\infty_T(L^2)}&\leq&\|b\|_{L^\infty_T(L^2)}\|b\|_{L^1_T(\mathcal X^{1})}\\
&\leq&2\varepsilon\|b^0\|_{L^2}\\
&\leq&\frac{\|b^0\|_{L^2}}{4}.
\end{array}$$
Then
\begin{equation}\label{eq401}
\|\psi(b)\|_{L^\infty_T(L^2)}\leq 2\|b^0\|_{L^2},\;\forall b\in B(\varepsilon,T).
\end{equation}
Finally, inequality (\ref{enq3}) gives
$$\begin{array}{lcl}
\|L(b)\|_{L^1_T(\mathcal X^1)}&\leq&\nu^{-1}\|a\|_{L^\infty_T(\mathcal X^{-1})}^{1/2}\|a\|_{L^1_T(\mathcal X^{1})}^{1/2}\|b\|_{L^\infty_T(\mathcal X^{-1})}^{1/2}\|b\|_{L^1_T(\mathcal X^{1})}^{1/2}\\
&\leq&\nu^{-1}\|a\|_{L^\infty_T(\mathcal X^{-1})}^{1/2}\|a\|_{L^1_T(\mathcal X^{1})}^{1/2}\sqrt{2}\sqrt{\varepsilon}\|b^0\|_{\mathcal X^{-1}}^{1/2}\\
&\leq&\varepsilon/3,\;(by\,(C3))\\\\
\|Q(b)\|_{L^1_T(\mathcal X^1)}&\leq&\nu^{-}\|b\|_{L^\infty_T(\mathcal X^{-1})}\|b\|_{L^1_T(\mathcal X^{1})}\\
&\leq&2\nu^{-1}\varepsilon\|b^0\|_{\mathcal X^{-1}}\\
&\leq&\varepsilon/3,\;(by\,(\ref{eq42})).
\end{array}$$
Then,
\begin{equation}\label{eq501}
\|\psi(b)\|_{L^1_T(\mathcal X^1)}\leq \varepsilon,\;\forall b\in B(\varepsilon,T).
\end{equation}
Therefore inequalities (\ref{eq301})-(\ref{eq401})-(\ref{eq501}) imply (\ref{eq45}).\\
\noindent\underline{\bf Proof of (\ref{eq46}):} Using inequality (\ref{enq1}), we obtain
$$\begin{array}{lcl}
\|L(\alpha_1)-L(\alpha_2)\|_{L^\infty_T(\mathcal X^{-1})}&\leq&\|L(\alpha_1-\alpha_2)\|_{L^\infty_T(\mathcal X^{-1})}\\
&\leq&\|a\|_{L^\infty_T(\mathcal X^{-1})}^{1/2}\|a\|_{L^1_T(\mathcal X^{1})}^{1/2}\|\alpha_1-\alpha_2\|_{L^\infty_T(\mathcal X^{-1})}^{1/2}\|\alpha_1-\alpha_2\|_{L^1_T(\mathcal X^{1})}^{1/2}\\
&\leq&\|a\|_{L^\infty_T(\mathcal X^{-1})}^{1/2}\|a\|_{L^1_T(\mathcal X^{1})}^{1/2}\|\alpha_1-\alpha_2\|_{\varepsilon,T}\\
&\leq&\frac{1}{12}\|\alpha_1-\alpha_2\|_{\varepsilon,T},\,(by\,(C8))\\\\
\|Q(\alpha_1)-Q(\alpha_2)\|_{L^\infty_T(\mathcal X^{-1})}&=&\|\int_0^te^{-\nu (t-z)\Delta}\mathbb P((\alpha_1-\alpha_2).\nabla \alpha_1+\alpha_2.\nabla (\alpha_1-\alpha_2))\|_{L^\infty_T(\mathcal X^{-1})}\\
&\leq&\displaystyle\Big(\sum_{i=1}^2\|\alpha_i\|_{L^\infty_T(\mathcal X^{-1})}^{1/2}\|\alpha_i\|_{L^1_T(\mathcal X^{1})}^{1/2}\Big)\|\alpha_1-\alpha_2\|_{L^\infty_T(\mathcal X^{-1})}^{1/2}\|\alpha_1-\alpha_2\|_{L^1_T(\mathcal X^{1})}^{1/2}\\
&\leq&2\sqrt{2}\sqrt{\varepsilon}\|b^0\|_{\mathcal X^{-1}}^{1/2}\|\alpha_1-\alpha_2\|_{\varepsilon,T}\\
&\leq&\frac{1}{12}\|\alpha_1-\alpha_2\|_{\varepsilon,T},\,(by\,(C7)).
\end{array}$$
Then
\begin{equation}\label{eq301}
\|\psi(\alpha_1)-\psi(\alpha_2)\|_{L^\infty_T(\mathcal X^{-1})}\leq \frac{1}{6}\|\alpha_1-\alpha_2\|_{\varepsilon,T},\;\forall \alpha_1,\alpha_2\in B(\varepsilon,T).
\end{equation}
Similarly, inequality (\ref{enq2}) gives
$$\begin{array}{lcl}
\|L(\alpha_1)-L(\alpha_2)\|_{L^\infty_T(L^2)}&\leq&\|L(\alpha_1-\alpha_2)\|_{L^\infty_T(L^2)}\\
&\leq&\|a\|_{L^\infty_T(L^2)}\|\alpha_1-\alpha_2\|_{L^1_T(\mathcal X^{1})}+\|a\|_{L^1_T(\mathcal X^{1})}\|\alpha_1-\alpha_2\|_{L^\infty_T(L^2)}\\
&\leq&\Big(\|a\|_{L^\infty_T(L^2)}+\|a\|_{L^1_T(\mathcal X^{1})}\Big)\|\alpha_1-\alpha_2\|_{\varepsilon,T}\\
&\leq&\Big(\|a^0\|_{L^2}+\|a\|_{L^1_T(\mathcal X^{1})}\Big)\|\alpha_1-\alpha_2\|_{\varepsilon,T}\\
&\leq&\Big(\|u^0\|_{L^2}+\|a\|_{L^1_T(\mathcal X^{1})}\Big)\|\alpha_1-\alpha_2\|_{\varepsilon,T}\\
&\leq&\Big(\frac{1}{24}+\|a\|_{L^1_T(\mathcal X^{1})}\Big)\|\alpha_1-\alpha_2\|_{\varepsilon,T},\;(by\,(\ref{L2}))\\
&\leq&\frac{1}{12}\|\alpha_1-\alpha_2\|_{\varepsilon,T},\;(by\,(C9))\\\\
\|Q(\alpha_1)-Q(\alpha_2)\|_{L^\infty_T(L^2)}&=&\|\int_0^te^{-\nu (t-z)\Delta}\mathbb P((\alpha_1-\alpha_2).\nabla \alpha_1+\alpha_2.\nabla (\alpha_1-\alpha_2))\|_{L^\infty_T(L^2)}\\
&\leq&\|\alpha_1-\alpha_2\|_{L^\infty_T(L^2)}\|\alpha_1\|_{L^1_T(\mathcal X^{1})}+\|\alpha_2\|_{L^\infty_T(L^2)}\|\alpha_1-\alpha_2\|_{L^1_T(\mathcal X^{1})}\\
&\leq&\Big(\|\alpha_1\|_{L^1_T(\mathcal X^{1})}+\|\alpha_2\|_{L^\infty_T(L^2)}\Big)\|\alpha_1-\alpha_2\|_{\varepsilon,T}\\
&\leq&(\varepsilon+2\|b^0\|_{L^2})\|\alpha_1-\alpha_2\|_{\varepsilon,T}\\
&\leq&(\varepsilon+2\|u^0\|_{L^2})\|\alpha_1-\alpha_2\|_{\varepsilon,T}\\
&\leq&\frac{1}{12}\|\alpha_1-\alpha_2\|_{\varepsilon,T},\;(by\,(\ref{L2})).
\end{array}$$
Then
\begin{equation}\label{eq701}
\|\psi(\alpha_1)-\psi(\alpha_2)\|_{L^\infty_T(L^2)}\leq \frac{1}{6}\|\alpha_1-\alpha_2\|_{\varepsilon,T},\;\forall \alpha_1,\alpha_2\in B(\varepsilon,T).
\end{equation}
Finally, inequality (\ref{enq3}) gives
$$\begin{array}{lcl}
\|L(\alpha_1)-L(\alpha_2)\|_{L^1_T(\mathcal X^1)}&=&\|L(\alpha_1-\alpha_2)\|_{L^1_T(\mathcal X^1)}\\
&\leq&\nu^{-1}\|a\|_{L^\infty_T(\mathcal X^{-1})}^{1/2}\|a\|_{L^1_T(\mathcal X^{1})}^{1/2}\|\alpha_1-\alpha_2\|_{L^\infty_T(\mathcal X^{-1})}^{1/2}\|\alpha_1-\alpha_2\|_{L^1_T(\mathcal X^{1})}^{1/2}\\
&\leq&\nu^{-1}\|a\|_{L^\infty_T(\mathcal X^{-1})}^{1/2}\|a\|_{L^1_T(\mathcal X^{1})}^{1/2}\|\alpha_1-\alpha_2\|_{\varepsilon,T}\\
&\leq&\frac{1}{12}\|\alpha_1-\alpha_2\|_{\varepsilon,T}\\\\
\|Q(\alpha_1)-Q(\alpha_2)\|_{L^1_T(\mathcal X^{1})}&=&\|\int_0^te^{-\nu (t-z)\Delta}\mathbb P((\alpha_1-\alpha_2).\nabla \alpha_1+\alpha_2.\nabla (\alpha_1-\alpha_2))\|_{L^1_T(\mathcal X^{1})}\\
&\leq&\displaystyle\nu^{-1}\Big(\sum_{i=1}^2\|\alpha_i\|_{L^\infty_T(\mathcal X^{-1})}^{1/2}\|\alpha_i\|_{L^1_T(\mathcal X^{1})}^{1/2}\Big)\|\alpha_1-\alpha_2\|_{L^\infty_T(\mathcal X^{-1})}^{1/2}\|\alpha_1-\alpha_2\|_{L^1_T(\mathcal X^{1})}^{1/2}\\
&\leq&\displaystyle\nu^{-1}\Big(\sum_{i=1}^2\|\alpha_i\|_{L^\infty_T(\mathcal X^{-1})}^{1/2}\|\alpha_i\|_{L^1_T(\mathcal X^{1})}^{1/2}\Big)\|\alpha_1-\alpha_2\|_{\varepsilon,T}\\
&\leq&\nu^{-1}2\sqrt{2}\sqrt{\varepsilon}\|b^0\|_{\mathcal X^{-1}}^{1/2}\|\alpha_1-\alpha_2\|_{\varepsilon,T}\\
&\leq&\frac{1}{12}\|\alpha_1-\alpha_2\|_{\varepsilon,T}.
\end{array}$$
Then,
\begin{equation}\label{eq801}
\|\psi(\alpha_1)-\psi(\alpha_2)\|_{L^1_T(\mathcal X^{1})}\leq \frac{1}{2}\|\alpha_1-\alpha_2\|_{\varepsilon,T},\;\forall \alpha_1,\alpha_2\in B(\varepsilon,T).
\end{equation}
Therefore inequalities (\ref{eq301})-(\ref{eq401})-(\ref{eq501}) gives
\begin{equation}\label{eq801}
\|\psi(\alpha_1)-\psi(\alpha_2)\|_{B(\varepsilon,T)}\leq \frac{1}{2}\|\alpha_1-\alpha_2\|_{B(\varepsilon,T)},\;\forall \alpha_1,\alpha_2\in B(\varepsilon,T).
\end{equation}
Fixed Point Theorem gives the existence and uniqueness of solution of $(RNS)$ in $C_T(L^2\cap \mathcal X^{-1})\cap L^1_T(\mathcal X^1).$ Therefore, we can  deduce the existence and  uniqueness of a local solution for Navier-Stokes system.
\section{\bf Proof of Theorem \ref{theo2}}
\noindent\underline{\bf Proof of (\ref{th2eq1}) :} In this subsection we want to prove the long time decay in $L^2\cap\mathcal X^{-1}$. Let $u\in C(\R^+,L^2\cap \mathcal X^{-1})$ be global solution of $(NS)$. By \cite{JB} we have
$$\limsup_{t\rightarrow\infty}\|u(t)\|_{\mathcal X^{-1}}=0.$$
Now, prove that $\limsup_{t\rightarrow\infty}\|u(t)\|_{L^2}=0.$ For a strictly positive real number $\delta$ and a given distribution $f$, we define the operators $A_\delta(D)$ and $B_\delta(D)$, respectively, by the following:
$$\begin{array}{lcl}
A_\delta(D)f&=&\mathcal F^{-1}\big({\bf 1}_{\{|\xi|<\delta\}}\widehat{f}\big),\\
B_\delta(D)f&=&\mathcal F^{-1}\big({\bf 1}_{\{|\xi|>\delta\}}\widehat{f}\big).
\end{array}$$
Let $u$ be a solution of $(NS)$. Denote by $ w_{\delta}=A_\delta(D)u$ and $ v_{\delta}=B_\delta(D)u$, respectively, the low-frequency part and the high-frequency part of $u$ and so on $w_{\delta}^0 $ and $ v_{\delta}^0$ for the initial data $u^0$. Applying the pseudo-differential operator $A_{\delta}(D)$ to the $(NS)$, we get
\begin{equation}\label{equ1}
\partial_{t}w_{\delta}-\nu\Delta w_{\delta} +A_\delta(D) \mathbb P(u.\nabla u)=0
\end{equation}
Taking the $L^2(\R^3)$-inner product and using the fact $A_\delta(D)^2=A_\delta(D)$, we obtain
$$\begin{array}{lcl}
\frac{1}{2}\frac{d}{dt}\|w_{\delta}(t)\|^2_{L^{2}}+\nu \|\nabla w_{\delta}(t)\|^2_{L^2}
&\leq&|\langle A_\delta(D) \mathbb P(u.\nabla u)(t)/w_{\delta}(t)\rangle_{L^{2}}|\\
&\leq&|\langle  \mathbb P(u.\nabla u)(t)/ A_\delta(D)w_{\delta}(t)\rangle_{L^{2}}|\\
&\leq&|\langle  \mathbb P(u.\nabla u)(t)/ w_{\delta}(t)\rangle_{L^{2}}|\\
&\leq&|\langle  u.\nabla u(t)/  \mathbb P( w_{\delta}(t))\rangle_{L^{2}}|\\
&\leq&|\langle  u.\nabla u(t)/   w_{\delta}(t)\rangle_{L^{2}}|\\
&\leq&|\langle ({\rm\,{div}}\,(u\otimes u))(t)/  w_{\delta}(t)\rangle_{L^{2}}|\\
&\leq&|\langle u\otimes u(t)/  \nabla w_{\delta}(t)\rangle_{L^{2}}|\\
&\leq&\|u\otimes u(t)\|_{L^{1}}\|\nabla w_{\delta}(t)\|_{L^\infty}\\
&\leq& (2\pi)^{-3}\|u(t)\|^{2}_{L^{2}}\|w_{\delta}(t)\|_{\mathcal X^{1}}
\end{array}$$
Integrating with respect to time and using Remark \ref{rem1}-(iv), we obtain
$$\|w_{\delta}(t)\|^2_{L^{2}}\leq \|w_{\delta}^0\|^2_{L^{2}}+  m_0\int_{0}^{t}\|w_{\delta}(s)\|_{\mathcal X^{1}}ds,$$
where $ m_0 =(2\pi)^{-3}\|u\|_{L^\infty(\R^+,L^2)}$. Also using Remark \ref{rem1}-(iii) we get $\|w_{\delta}(t)\|^2_{L^{2}}\leq M_\delta, $ where
$$M_\delta = \|w_{\delta}^0\|^2_{L^{2}}+ m_0\int_{0}^{\infty}\|w_{\delta}(s)\|_{\mathcal X^{1}}ds .$$
On the one hand, it is clear that $\lim_{\delta\rightarrow 0}\|w_{\delta}^0\|^2_{L^{2}}=0.$ On the other, we have  $\lim_{\delta\rightarrow 0} \| w_{\delta(t)}\|_{\mathcal X^{1}}=0$ and  $\| w_{\delta(t)}\|_{\mathcal X^{1}}\leq\|u(t)\|_{\mathcal X^{1}} \in  L^{1}([0,\infty)).$ Then Dominate Convergence Theorem implies that
$$\lim_{\delta\rightarrow 0}\int_{0}^{\infty}\|w_{\delta}(s)\|_{\mathcal X^{1}}ds = 0.$$
Hence,  $\lim_{\delta\rightarrow 0}M_\delta =0$ and thus
\begin{equation}\label{uq1}
\lim_{\delta\rightarrow 0}\sup_{t\geq0}\|w_\delta(t)\|_{L^2}^{2}\rightarrow0.
\end{equation}
Let us investigate the high-frequency part. To do so, one applies the pseudo-differential operator $B_{\delta}(D)$ to the $(NS)$ to get
\begin{equation}\label{equ4}
\partial_{t}v_{\delta}-\nu\Delta v_{\delta} +B_\delta(D) \mathbb P(u.\nabla u)=0.
\end{equation}
The integral form of $v_{\delta}$ is
$$v_{\delta}(t) = e^{\nu t\Delta}v_{\delta}^0-\int_0^t e^{\nu(t-\tau)\Delta}B_\delta(D) \mathbb P(u.\nabla u)d\tau.$$
Taking the $L^2(\R^3)$ norm, we obtain
$$\begin{array}{lcl}
\|v_{\delta}(t)\|_{L^2}&\leq& \|e^{\nu t\Delta}v_{\delta}^0\|_{L^2}+\int_0^t\| e^{\nu(t-\tau)\Delta}B_\delta(D)\mathbb P(u.\nabla u)\|_{L^2}d\tau \\
&\leq& e^{-\nu t\delta^2}\|v_{\delta}^0\|_{L^2}+\int_0^t e^{-\nu(t-\tau)\delta^2}\|u\nabla u\|_{L^2}d\tau \\
&\leq& e^{-\nu t\delta^2}\|u^0\|_{L^2}+\int_0^t e^{-\nu(t-\tau)\delta^2}\|u\|_{L^2}\|\nabla u\|_{L^\infty}d\tau .
\end{array}$$
Then $$\|v_{\delta}(t)\|_{L^2} \leq e^{-\nu t\delta^2}\|u^0\|_{L^2}+m_{0}\int_0^t e^{-\nu(t-\tau)\delta^2}\|u(\tau)\|_{\mathcal X^{1}} d\tau:= G_{\delta}(t).$$
We have
$$\int_{0}^{\infty}G_{\delta}(t)dt\leq\frac{\|u^0\|_{L^2}}{\nu \delta^2}+\frac{m_{0}}{\nu \delta^2}\|u\|_{L^1(\R^+,\mathcal{X}^1 )}<\infty.$$
This leads to the fact that the function  $(t \rightarrow \|v_{\delta}(t)\|_{L^2})$ is both continuous and Lebesgue integrable over $\R^+$. Let $\varepsilon>0$ be positive real number small enough. Firstly, equation (\ref{uq1} ) implies that some $\delta _{\varepsilon} > 0$ exists such that
\begin{equation}\label{equ5}
 \|w_{\delta_{\varepsilon}}(t)\|_{L^2} \leq\varepsilon/2, \;\;\;\forall t\geq0 .
 \end{equation}
Secondly, consider the set $R_{\delta_{\varepsilon}}$ defined by
\begin{equation}\label{eqRd1}R_{\delta_{\varepsilon}} :=\{t>0, \|v_{\delta_{\varepsilon}}(t)\|_{L^2} >\varepsilon/2 \}.\end{equation}  If we denote by $ \lambda_{1}(R_{\delta_\varepsilon})$ the Lebesgue measure of $R_{\delta_\varepsilon} $, we have
$$\int_{0}^{\infty} \|v_{\delta_{\varepsilon}}(t)\|_{L^2}dt\geq\int_{R_{\delta_\varepsilon} } \|v_{\delta_\varepsilon}(t)\|_{L^2{(\R^3)}}dt\geq \frac{\varepsilon}{2}\lambda_{1}(R_{\delta_\varepsilon}). $$
 By this, we can deduce that $\lambda_{1}(R_{\delta_\varepsilon})\leq T_{\varepsilon},$  where $ T_{\varepsilon} =(2/\varepsilon)\int_{0}^{\infty} \|v_{\delta_{0}}(t)\|_{L^2{(\R^3)}}dt.$ Then, there is $t_\varepsilon \in[0,T_{\varepsilon} + 1]$ such that $t_{\varepsilon}$ does not belong to $R_{\delta_{\varepsilon}}.$ This implies that
 \begin{equation}\label{equ6}
 \|v_{\delta_{\varepsilon}}(t_\varepsilon)\|_{L^2{(\R^3)}}\leq\varepsilon/2.
 \end{equation}
Equations (\ref{equ5}) and (\ref{equ6}) together with triangular inequality imply that $ \|u(t_\varepsilon)\|_{L^2(\R^3)}< \varepsilon.$
For $t\geq t_\varepsilon$, we have
$$\begin{array}{lcl}
\|u(t)\|_{L^2}&\leq&\|u(t_\varepsilon)\|_{L^2}\exp((2\pi)^{-3}\int_{t_0}^\infty \|u(z)\|_{\mathcal X^{1}}dz)  \\
&\leq&\varepsilon\exp((2\pi)^{-3}\int_0^\infty \|u\|_{\mathcal X^{1}}) .
 \end{array}$$
It suffices to replace $\varepsilon $ by $\varepsilon\exp(-(2\pi)^{-3}\int_0^\infty \|u\|_{\mathcal X^{1}})$ in (\ref{equ5})-(\ref{eqRd1})-(\ref{equ6}) we get the desired result.\\
\noindent\underline{\bf Proof of (\ref{th2eq2}) :} In this subsection we want to give a precision for the decay of $\|u(t)\|_{\mathcal X^{-1}}$ at $\infty$.
Let $\varepsilon>0$ such that $\varepsilon<\epsilon_0$($\epsilon_0$ is given by Theorem \ref{th14}), by (\ref{th2eq1})  we can  suppose that,
$$\|u^0\|_{\mathcal X^{-1}}  <\min(\varepsilon,\frac{\nu}{2})\;\;{\rm{and}}\;\;\|u^0\|_{L^2} < \varepsilon/2.$$
Then, by Remark\ref{rem1}-(ii)-(iv) we get $\|u(t)\|_{L^2}\leq2\|u^0\|_{L^2}< \varepsilon$ for all $t\geq0$ and
\begin{equation}\label{eqx1}\|u(t)\|_{\mathcal X^{-1}}+\frac{\nu}{2}\int_0^t\|u(z)\|_{\mathcal X^{1}}dz\leq\|u^0\|_{\mathcal X^{-1}}<\frac{\nu}{2},\;\;\forall t\geq0.\end{equation}
For $\lambda >0$ and $t>00$, we have
$$\|u(t)\|_{\mathcal X^{-1}}=I_\lambda(t)+J_\lambda(t),$$
with
$$\begin{array}{lcl}
I_\lambda(t)&=&\int_{\{\xi\in\R^3/|\xi|<\lambda\}}\frac{|\widehat{u}(t,\xi)|}{|\xi|}d\xi\\
J_\lambda(t)&=&\int_{\{\xi\in\R^3/|\xi|>\lambda\}}\frac{|\widehat{u}(t,\xi)|}{|\xi|}d\xi.
\end{array}$$
We have
$$\begin{array}{lcl}
I_\lambda(t)&\leq&\Big(\int_{\{\xi\in\R^3/|\xi|<\lambda\}}\frac{1}{|\xi|^2}d\xi\Big)^{1/2}\|\widehat{u}\|_{L^2}\\
&\leq& c_{0}\Big(\int_{0}^{\lambda} dr \Big)^{1/2}\|\widehat{u}\|_{L^2}\\
&\leq& c_{0} \sqrt{\lambda}\| \widehat{u}(t)\|_{L^2}\\
&\leq& c_{0} \sqrt{\lambda}\| \widehat{u^0}\|_{L^2}\\
&\leq& c_{1} \sqrt{\lambda}\|u^0\|_{L^2}
\end{array}$$
and
$$\begin{array}{lcl}
J_\lambda(t)&\leq&\int_{\{\xi\in\R^3/|\xi|>\lambda\}}e^{-\nu\sqrt{t/2}|\xi|}e^{\sqrt{t/2}|\xi|}\frac{|\widehat{u}(t,\xi)|}{|\xi|}d\xi\\
&\leq&  e^{-\sqrt{\nu t/2}\lambda}  \int_{\R^3}  e^{\sqrt{\nu t/2}|\xi|}\frac{|\widehat{u}(t,\xi)|}{|\xi|}d\xi.
\end{array}$$
For a fixed time $t>0$ the  $ v:(z,x)\rightarrow u(\frac{t}{2}+z,x)$ satisfies $\|v(0)\|_{\mathcal X^{-1}}<\epsilon_0$ and it is the unique global solution of the following system,
 $$
  \begin{cases}
     \partial_t v
 -\nu\Delta v+ v.\nabla v   =\;\;-\nabla q \\
    v(0,x) =u(\frac{t}{2},x).
  \end{cases}
$$
By Theorem \ref{th14}, we get
$$\int_{\R^3}  e^{\sqrt{\nu z}|\xi|}\frac{|\widehat{v}(z,\xi)|}{|\xi|}d\xi +\frac{\nu}{2} \int_{0}^{z} \int_{\R^3}  e^{\sqrt{\nu s}|\xi|}\frac{|\widehat{v}(s,\xi)|}{|\xi|}d\xi ds \leq 2\|v(0)\|_{\mathcal X^{-1}}$$
or
$$\int_{\R^3}  e^{\sqrt{\nu z}|\xi|}\frac{|\widehat{u}(\frac{t}{2}+z,\xi)|}{|\xi|}d\xi +\frac{\nu}{2} \int_{0}^{z} \int_{\R^3}  e^{\sqrt{\nu s}|\xi|}\frac{|\widehat{u}(\frac{t}{2}+s,\xi)|}{|\xi|}d\xi ds \leq 2\int_{\R^3} \frac{|\widehat{u}(\frac{t}{2},\xi)|}{|\xi|}d\xi.$$
For $ z = \frac{t}{2}$,  we get $ \int_{\xi}  e^{\sqrt{\nu t/2}|\xi|}\frac{|\widehat{u}(t,\xi)|}{|\xi|}d\xi \leq \|u(t/2)\|_{\mathcal X^{-1}}$, which implies
$$J_\lambda(t) \leq e^{-\sqrt{\nu t/2}\lambda }\|u(t/2)\|_{\mathcal X^{-1}}.$$
Then
$$\|u(t)\|_{\mathcal X^{-1}} \leq c_1 \sqrt{\lambda}\| u^0\|_{L^2}+e^{-\sqrt{\nu t/2}\lambda }\|u(t/2)\|_{\mathcal X^{-1}}.$$
Multiplying this inequality by $t^{1/4}$
$$ t^{1/4}\|u(t)\|_{\mathcal X^{-1}} \leq  t^{1/4} c_1 \sqrt{\lambda}\| u^0\|_{L^2}+  2^{1/4}e^{-\sqrt{\nu t/2}\lambda } (\frac{t}{2})^{1/4}\|u(t/2)\|_{\mathcal X^{-1}}$$
and choose $\lambda>0$ such that $$2^{1/4}e^{-\sqrt{\nu t/2}\lambda }=1/2  \Rightarrow \sqrt {\nu t/2}\lambda =5/4 \ln{2}\Rightarrow \lambda = \frac{5\sqrt{2}\ln{2}}{4\sqrt{\nu t}}$$
we obtain
$$t^{1/4}\|u(t)\|_{\mathcal X^{-1}} \leq M_{0} +\frac{1}{2}  (\frac{t}{2})^{1/4} \|u(t/2)\|_{\mathcal X^{-1}}$$ with $$M_{0}= c_{0}(\frac{5\sqrt{2}\ln{2}}{4\sqrt{\nu}})^{1/2}\| u^0\|_{L^2}.$$
Applying Lemma \ref{lemt1} and Remark \ref{remlemt1} with
$$f(t)= t^{1/4}\|u(t)\|_{\mathcal X^{-1}}, \;\;\;\; \theta_{1}=\theta_{2}=1/2, $$
we get $$\limsup_{t\rightarrow +\infty}  t^{1/4}\|u(t)\|_{\mathcal X^{-1}}\leq  2 M_{0}. $$
Applying this result to the solution of the following system, for $a\geq0 $
$$
  \begin{cases}
\partial_t w
 -\nu\Delta w+ w.\nabla w    =\;\;-\nabla h \hbox{ in } \mathbb R^+\times \mathbb R^3\\
 {\rm div}\, v = 0 \hbox{ in } \mathbb R^+\times \mathbb R^3\\
   w(0,x) = u(a,x)  \hbox{ in } \mathbb R^3,
  \end{cases}
$$
we obtain $$\limsup_{t\rightarrow +\infty}  t^{1/4}\|u(t)\|_{\mathcal X^{-1}}\leq c_{0}(\frac{5\sqrt{2}\ln{2}}{4\sqrt{\nu}})^{1/2}\| u(a)\|_{L^2}.$$
Then the fact $\lim_{a\rightarrow +\infty} \| u(a)\|_{L^2}=0$  implies the desired result.
\section{\bf Long time decay in $\mathcal X^\sigma$} In this section we want to prove Corollary \ref{cor1}.\\
\noindent\underline{First case}\;\; $:-3/2<\sigma<-1.$ For $\lambda >0$ and $t>0$, we have
$$\|u(t)\|_{\mathcal X^{\sigma}}=I_{1}(t,\lambda)+I_{2}(t,\lambda)$$
$$\begin{array}{lcl}
I_{1}(t,\lambda)&=&\int_{\{\xi\in\R^3/|\xi|<\lambda\}}|\xi|^\sigma|\widehat{u}(t,\xi)|d\xi\\
I_{2}(t,\lambda)&=&\int_{\{\xi\in\R^3/|\xi|>\lambda\}}|\xi|^\sigma|\widehat{u}(t,\xi)|d\xi.
\end{array}$$
We have
$$\begin{array}{lcl}
I_{1}(t,\lambda)&\leq& \Big(\int_{|\xi|<\lambda}|\xi|^{2\sigma} d\xi \Big)^{1/2}\|\widehat{u}(t)\|_{L^2}\\
&\leq& c_1 \lambda^{\sigma+3/2}\| u(t)\|_{L^2}
\end{array}$$
and
$$\begin{array}{lcl}
I_{2}(t,\lambda)&\leq&\int_{\{\xi\in\R^3/|\xi|>\lambda\}}|\xi|^{\sigma+1}\frac{|\widehat{u}(t,\xi)|}{|\xi|}d\xi\\
&\leq& \lambda^{\sigma+1}\int_{\R^3}\frac{|\widehat{u}(t,\xi)|}{|\xi|}d\xi\\
&\leq& \lambda^{\sigma+1}\|u(t)\|_{\mathcal{ X}^{-1}} .
\end{array}$$
We get $\|u(t)\|_{\mathcal X^{\sigma}}=A\lambda^{\sigma+3/2}+B\lambda^{\sigma+1}:=\varphi(\lambda)$, with
$$ A = c_{0}\|{u}(t)\|_{L^2}\;\;{\rm and}\;\;B=\|u(t)\|_{\mathcal{ X}^{-1}}.$$
The study of the function $\varphi$ gives
$$\varphi'(\lambda)=(\sigma+3/2) A \lambda^{\sigma+1/2}+(\sigma+1) B \lambda^{\sigma},$$
then$$\varphi'(\lambda)=0\Rightarrow \lambda=\lambda_{0}=(\frac{-(1+\sigma)B}{(\sigma+3/2)A})^{2}.$$
For $ \lambda=\lambda_{0}$, we get
$$\begin{array}{lcl}
\|u(t)\|_{\mathcal X^{\sigma}}&\leq&A(\frac{-(1+\sigma)B}{
(\sigma+3/2)A})^{2\sigma+3}+B(\frac{-(1+\sigma)B}{(\sigma+3/2)A})^{2\sigma+2}\\
&\leq&c_{\sigma}A^{-2\sigma-2}B^{3+2\sigma}.
\end{array}$$
Then
\begin{equation*}
\|u(t)\|_{\mathcal X^{\sigma}} \leq c'_{\sigma}(\|{u}(t)\|_{L^2})^{-2\sigma-2}(\|u(t)\|_{\mathcal{ X}^{-1}})^{3+2\sigma}.
\end{equation*}
Theorem \ref{th15} implies$$\|u(t)\|_{\mathcal X^{-1}}=o(t^{-1/4})\;\;{\rm and}\;\;\|{u}(t)\|_{L^2}\rightarrow 0,$$
which gives the desired result.\\
\noindent\underline{Second case} \;\; $:-1<\sigma.$ By Theorem \ref{th15} we can assume that $\|u^0\|_{\mathcal X^{-1}}<\epsilon_0$ and Theorem \ref{th14} gives,
$$\begin{array}{lcl}
\|u(t)\|_{\mathcal X^{\sigma}}&=&\int_{\R^3} e^{-\sqrt{\nu t/2}|\xi|}|\xi|^{\sigma+1}e^{\sqrt{\nu t/2}|\xi|}\frac{|\widehat{u}(t,\xi)|}{|\xi|}d\xi\\
&=&\frac{1}{(\sqrt{t/2})^{\sigma+1}}\int(\sqrt{\nu t/2}|\xi|)^{\sigma+1} e^{-\sqrt{\nu t/2}|\xi|} e^{\sqrt{\nu t/2}|\xi|}\frac{|\widehat{u}(t,\xi)|}{|\xi|}d\xi\\
&\leq& C_\nu t^{-\frac{\sigma+1}{2}}\int e^{\sqrt{\nu t/2}|\xi|}\frac{|\widehat{u}(t,\xi)|}{|\xi|}d\xi\\
&\leq& 2C_\nu t^{-\frac{\sigma+1}{2}}\|u(t/2)\|_{\mathcal X^{-1}},
\end{array}$$
with $ C_\nu=\nu^{-\frac{\sigma+1}{2}}\sup_{z\geq0}z^{\sigma+1} e^{-rz}.$ Combining this result with the fact $\|u(t/2)\|_{\mathcal X^{-1}}=o(t^{-1/4}) $
we get the desired result.

\end{document}